\documentclass[11pt,reqno]{amsart}
\usepackage{amsthm,amsfonts,amssymb,euscript}

\usepackage{hyperref}
\usepackage{xcolor}
\hypersetup{
colorlinks=true,
linkcolor=blue,    
citecolor=red,   
urlcolor=cyan      
}

\usepackage{comment}

\DeclareMathOperator*{\med}{med}

\def\normo#1{\left\|#1\right\|}

\def\norm#1{\|#1\|}

\newcommand{\X}{{\widehat{X}}}
\newcommand{\N}{{\mathbb N}}

\newcommand{\R}{{\mathbb R}}

\newcommand{\Z}{{\mathbb Z}}
\newcommand{\ft}{{\mathcal{F}}}

\newcommand{\les}{{\lesssim}}
\newcommand{\ges}{{\gtrsim}}
\newcommand{\ra}{{\rightarrow}}

\newcommand{\Sch}{{\mathcal{S}}}
\newcommand{\supp}{{\mbox{supp}}}

\theoremstyle{plain}
\newtheorem{theorem}[subsection]{Theorem}
\newtheorem{proposition}[subsection]{Proposition}
\newtheorem{lemma}[subsection]{Lemma}

\theoremstyle{remark}
\newtheorem{remark}[]{Remark}

\theoremstyle{definition}

\numberwithin{equation}{section}

\begin{document}
\title[LWP of dgBO in $\widehat{H}^{s}_{r}(\mathbb{R})$]{Local well-posedness for dispersion generalized Benjamin-Ono equations in Fourier-Lebesgue spaces}

\author[Z. Chen]{Zijun Chen}
\address{School of Mathematics, Monash University, VIC 3800, Australia }
\email{zijun.chen@monash.edu}


	
\begin{abstract}
We prove that the Cauchy problem for the dispersion
generalized Benjamin-Ono equation where $0<\alpha \leq 1$
\begin{eqnarray*}
\left\{
\begin{array}{l}
\partial_t u+|\partial_x|^{1+\alpha}\partial_x u+uu_x=0,\\
u(x,0)=u_0(x),
\end{array}
\right.
\end{eqnarray*}
is locally well-posed in the Fourier-Lebesgue space $\widehat{H}^{s}_{r}(\mathbb{R})$. This is proved via Picard iteration arguments using $X^{s,b}$-type space adapted to the Fourier-Lebesgue space, inspired by the work of Gr\"unrock and Vega. Note that, previously, Molinet, Saut and
Tzvetkov \cite{MST2001} proved that the solution map is not $C^2$ in $H^s$ for any $s$ if $0\leq \alpha<1$. However, in the Fourier-Lebesgue space, we have a stronger smoothing effect to handle the $high\times low$ interactions.
\end{abstract}

\maketitle
	
\tableofcontents
	
\section{Introduction}
	
In this paper, we consider the Cauchy problem for the dispersion
generalized Benjamin-Ono equation (dgBO)
\begin{eqnarray}\label{eq:dgBO}
\left\{
\begin{array}{l}
\partial_t u+|\partial_x|^{1+\alpha}\partial_x u+uu_x=0,\ (x,t)\in \R^2,\\
u(x,0)=u_0(x),
\end{array}
\right.
\end{eqnarray}
where $0\leq \alpha \leq 1$. Here $u: \R^2 \ra \R$ is a real-valued
function and $|\partial_x|$ is the Fourier multiplier operator with
symbol $|\xi|$. These equations arise as mathematical models for the weakly nonlinear propagation of long waves. The case $\alpha=0$ corresponds to the Benjamin-Ono equation and the case $\alpha=1$ corresponds to the Korteweg-de Vries equation. 
The equation \eqref{eq:dgBO} is invariant under the following scaling transform:
\begin{equation}\label{scaling}
u(x, t) \rightarrow u_\lambda(x, t)=\lambda^{\frac{1+\alpha}{2}} u\left(\lambda x, \lambda^{2+\alpha} t\right), \quad \lambda>0.
\end{equation}
The scaling critical Sobolev space of \eqref{eq:dgBO} is $\dot{H}^{-\frac{1+\alpha}{2}}$ in the sense that the homogeneous Sobolev norm is invariant under the scaling transform \eqref{scaling}. 

The low regularity well-posedness of \eqref{eq:dgBO} in Sobolev spaces $H^{s}(\mathbb{R})$ has been extensively studied in recent years; we will summarize the most recent findings.

$\alpha=0$: \eqref{eq:dgBO} is known as the Benjamin-Ono equation. Ponce \cite{P1991} first proved $C^{0}$ local well-posedness ($C^{0}$ LWP) in $H^{s}$ for $s\geq \frac{3}{2}$ using the energy method with enhanced dispersive smoothing. Later, Koch-Tzvetkov \cite{KT2003} refined the energy method and smoothing effect to improve LWP for $s>\frac{5}{4}$. This approach was further refined by Kenig-Keonig \cite{KK2003}, which led to LWP for $s>\frac{9}{8}$.

Tao \cite{T2004} applied a gauge transformation to effectively remove the derivative (or at least the worst terms that involve the derivative) from the nonlinearity and obtained the global well-posedness (GWP) for $s\geq 1$. By a synthesis of Tao’s gauge transformation and $X^{s, b}$ techniques, Ionescu-Kenig \cite{IK2007} improved the range to $s\geq 0$. Recently, Killip-Laurens-Visan \cite{KLV2023} obtained a sharp GWP for $s>-\frac{1}{2}$ by using a new gauge transform and delving into the complete integrability of BO equations.

$\alpha=1$: \eqref{eq:dgBO} is the KdV equation. The first $C^{2}$ local well-posedness by contraction principle was proved by Kenig-Ponce-Vega \cite{KPV1993} for $s>\frac{3}{4}$. Bourgain \cite{B1993} extended this result to $s \geq 0$ by developing $X^{s, b}$ spaces. Then by developing the bilinear estimates in $X^{s, b}$ space, Kenig-Ponce-Vega \cite{KPV1996} were able to prove local well-posedness for $s>-\frac{3}{4}$. 
Christ-Colliander-Tao \cite{CCT2003}, Guo \cite{G2009} and Kishimoto \cite{K2009} discussed the endpoint $s=-\frac{3}{4}$. Recently, Killip-Visan \cite{KV2019} obtained a sharp $C^{0}$ GWP for $s\geq-1$ by using similar strategies mentioned above for BO equations.

$0<\alpha<1$:
Kenig-Ponce-Vega \cite{KPV1991} have shown that \eqref{eq:dgBO} is $C^{0}$ locally well-posed provided $s \geq \frac{3}{4}(2-\alpha)$, using the energy method with enhanced smoothing effect.
The Sobolev index has been pushed down to $s>1-\alpha$ by Guo \cite{G2012}.
Herr-Ionescu-Kenig-Koch \cite{HIKK2010} used a para-differential renormalization method to show the range of $s\geq 0$.

Note that there is only $C^{0}$ local well-posedness in $H^s$ if $0 \leq \alpha<1$. This indicates that only the $H^s$ assumption on the initial data is not sufficient to prove local well-posedness of dgBO via Picard iteration, as the ill-posedness result from Molinet-Saut-Tzvetkov \cite{MST2001}, by showing the solution mapping fails to be $C^2$ smooth from $H^s$ to $C\left([0, T]; H^s\right)$ at the origin for any $s$.
The reason is that the dispersive effect of the dispersive group of dgBO is too weak to recover the derivative in the nonlinearity. Hence the {$high\times low$} interactions break down the $C^2$ smoothness. 

Therefore, to study the well-posedness of dgBO when $0 \leq \alpha \leq 1$ in the sense that the solution mapping is uniformly continuous, we might choose to abandon $H^s$ and prove it still via contraction mapping principle in some other space of initial data. 
For instance, Herr \cite{H2007} applied a weighted Sobolev data space to obtain the well-posedness of dgBO equations for $0<\alpha<1$.
We also got some inspiration from the work of Gr\"unrock and Vega. Gr\"unrock \cite{G2004} obtained LWP of the modified KdV equation in Fourier-Lebesgue spaces $\widehat{H}^{s}_{r}$ for $\frac{4}{3} <r\leq 2$ and $s \geq \frac{1}{2}-\frac{1}{2r}$, which was enhanced by Gr\"unrock-Vega \cite{GV2009} to $1 <r\leq 2$ and $s \geq \frac{1}{2}-\frac{1}{2r}$. 
Furthermore, Gr\"unrock in \cite{Gr2010} studied the hierarchies of higher order mKdV and KdV equations systematically by using a mixed resolution space $\X^{s, b}_{r, p}$ where the time parameter $p$ depends on the value of $r$. 
This initial data space, Fourier-Lebesgue space $\widehat{H}^{s}_{r}$, involves $L^r$-type integrability to a spatial weight, in addition to $H^s$ regularity, where the norm is defined by
\begin{equation}\label{spa_FL}
\|f\|_{\widehat{H}_{r}^{s}(\mathbb{R})}:=\left\|\langle\xi\rangle^{s} \widehat{f}(\xi)\right\|_{L_{\xi}^{r^{\prime}}}, 1/r+1/r^{\prime}=1
\end{equation}
for $s \in \mathbb{R}, 1 \leq r \leq \infty$. From the scaling point of view, the spaces $\widehat{H}^{s}_{r}$ behave like the Bessel potential spaces ${H}^{s}_{r}$ which are embedded in $\widehat{H}^{s}_{r}$ for $1 \leq r \leq 2$ by Hausdorff-Young inequality, and like $H^\sigma$ if $s-\frac{1}{r}+\frac{1}{2}=\sigma$. 

Now, we state our main result:
\begin{theorem}\label{theorem}
Let $0 <\alpha \leq 1$ and $1 < r < 1+\alpha$.
The Cauchy problem \eqref{eq:dgBO} is locally well-posed in $\widehat{H}_{r}^{s}(\mathbb{R})$ if 
\begin{equation}\label{regularity}
s > -1-\alpha+\frac{2}{r}+\frac{\alpha}{2r}.
\end{equation}
\end{theorem}
This paper emphasises the local well-posedness result, which can be obtained by the
contraction mapping principle so that the flow map is real analytic. Since the dispersive effect of dgBO is weak when $0\leq \alpha <1$, it previously failed to prove the local well-posedness by Picard iteration in the classical Sobolev space $H^s$. However, we observe a stronger local smoothing effect in $\widehat{H}_{r}^{s}$ when $r<2$, hence we may expect something better in this type of space. The answer is yes and the contraction mapping principle works.

\section{Function spaces and linear estimates}
We use the notation $X\lesssim Y$ for $X, Y \in \mathbb{R}$ to denote that there exists a constant $C>0$ such that $X\leq CY$.  The notation $X\sim Y$ denotes that there exist positive constants $c, C$ such that $cY \leq X \leq C Y$. 
For $a\in \R$, $a\pm$ denotes $a\pm \varepsilon$ for any sufficiently small $\varepsilon>0$.
We use capitalized variables $\{N, L, N_1, N_2,\cdots\}$ to denote dyadic numbers, unless otherwise specified.
$\widehat{u}$ denotes the standard Fourier transform $\mathcal{F}_x u$, $\mathcal{F}_t u$ or $\mathcal{F}_{t, x} u$. 
Let $\omega(\xi)=-\xi|\xi|^{1+\alpha}$ be the dispersion relation associated with the equation \eqref{eq:dgBO} and $W_{\alpha}(t)=\ft^{-1}e^{it\omega(\xi)}\ft$ be the linear propagator.

Let $\psi \in C_0^\infty(\R)$ be a real-valued, non-negative, even, and radially-decreasing function such that $\supp \psi\subset [-5/4, 5/4]$ and $\psi\equiv 1$ in $[-1, 1]$. 
$\Z_+=\N \cup \{0\}$.
For a dyadic number $N\in 2^{\Z_+}$, denote $\chi(\xi):=\psi(\xi)-\psi(2\xi)$ and $\chi_{N}:=\chi(N^{-1}\cdot)$. The Littlewood-Paley projectors for frequency and modulation are defined by
\begin{equation}\label{dyadic}
\begin{aligned}
\widehat{P_1f}=\psi(\xi)\widehat{f}, &\quad \widehat{P_Nf}=\chi_N(\xi)\widehat{f}\ \text{for} \ N\geq 2;\\
\widehat{Q_1f}=\psi(\tau-\omega(\xi))\widehat{f}, &\quad \widehat{Q_Lf}=\chi_L(\tau-\omega(\xi))\widehat{f} \ \text{for} \ L\geq 2.
\end{aligned}
\end{equation}
The Fourier-Lebesgue type Bourgain space $\widehat{X}_r^{s,b}$ associated to \eqref{eq:dgBO} is defined by the norm
\begin{equation}\label{Fourier restriction norm}
\|u\|_{\widehat{X}_r^{s,b}}=\left\|\langle \xi \rangle^s \langle \tau-\omega(\xi)\rangle^b \widehat{u}(\xi, \tau)\right\|_{L^{r^{\prime}}_{\tau, \xi}}
\end{equation}
for $s, b \in \mathbb{R}$ and $1 \leq r \leq \infty$, where $\langle \cdot \rangle=(1+|\cdot|^2)^{\frac{1}{2}}$. When $s=b=0$, we write $\widehat{X}^{s,b}_{r}$ as $\widehat{L}^{r}_{t, x}$ for simplicity. We use the dyadic frequency localization operators $P_N$ and $Q_L$ to rewrite \eqref{Fourier restriction norm} as
\begin{equation}
\|u\|_{\widehat{X}^{s, b}_{r}}=\left(\sum_{N, L} N^{s r^\prime } L^{b r^\prime }\left\|P_N Q_L u\right\|_{\widehat{L}_{t, x }^r}^{r^\prime}\right)^{\frac{1}{r^\prime}}.
\end{equation}

By slightly modifying the proof of mKdV equations for the unitary group $\{e^{-t\partial_x^{3}}\}$ in \cite{G2004}, it is easy to obtain some linear estimates for the unitary group $\{e^{-t|\partial_x|^{1+\alpha}}\}$ in $\widehat{X}^{s, b}_{r}$ spaces.
\begin{lemma}[Extension Lemma]\label{lem:ext}
Let  $Z$ be any space-time Banach space satisfying the time modulation estimate
\begin{equation}
\left\|g(t)F(t,x)\right\|_Z\leq \|g\|_{L_t^{\infty}(\R)}\norm{F(t,x)}_Z 
\end{equation}
for any $F\in Z$ and $g\in L_t^{\infty}(\R)$. Let $T:(h_1, \cdots, h_k)\rightarrow T(h_1, \cdots, h_k)$ be a spatial multilinear operator for which one has the estimate
\begin{equation}
\left\|T(W_{\alpha}(t)f_{1},\cdots, W_{\alpha}(t)f_{k})\right\|_Z\les \prod_{j=1}^k \|f_{j}\|_{\widehat{L}_x^r}
\end{equation}
for all $f_{1},\cdots, f_{k}\in \widehat{L}^r_x$. Then for $b>1/r$, we have the estimate
\begin{equation}
\|T(u_1, \cdots, u_k)\|_{Z}\les_k \prod_{j=1}^k\|u_j\|_{\widehat{X}_r^{0,b}}
\end{equation}
for all $u_1, \cdots, u_k \in \widehat{X}_r^{0,b}$.
\end{lemma}
\begin{proof}
When $k=1$, this is proved in Lemma 2.1 \cite{G2004}. For $k\geq 2$, one can prove it by slightly modifying the proof in \cite{T2007} Lemma 4.1 and \cite{G2004} Lemma 2.1. We omit the details. 
\end{proof}
As a corollary of Lemma \ref{lem:ext}, if $b>1/r$, we have the following embedding
\begin{equation}
\widehat{X}_{r}^{s, b}\subset C(\mathbb{R}; \widehat{H}^{s}_{r}).
\end{equation}

\begin{lemma}[Strichartz estimates]\label{Strichartz}
Assume that $(q, p)$ satisfies 
\begin{equation}
\frac{2}{q}+\frac{1}{p}=\frac{1}{r}.
\end{equation}
Then for all $\phi \in \Sch(\R)$
\begin{eqnarray}
\norm{D^{\alpha/q} W_\alpha(t)\phi}_{L_t^qL_x^p}\les \|\phi\|_{\widehat{L}_{x}^{r}}
\end{eqnarray}
holds if one of the following conditions is satisfied:
\begin{equation}
\begin{aligned}
& (1)\ 4 \leq q \leq \infty, 4 < p \leq \infty; \\
& (2)\ 1 / 4 \leq 1 / p \leq 1 / p+1 / q<1 / 2; \\
& (3)\ (q, p)=(\infty, 2).
\end{aligned}
\end{equation}
Moreover, by Lemma \ref{lem:ext}, we have for any $b>\frac{1}{r}$ and $u \in \widehat{X}_r^{0, b}$,
\begin{equation}
\left\|D^{\alpha / q} u\right\|_{L_t^q L_x^p} \lesssim\|u\|_{\widehat{X}_r^{0, b}}.
\end{equation}
\end{lemma}

\begin{lemma}[Local smoothing estimates]\label{localsmooth}
For all $\phi \in \Sch(\R)$, we have
\begin{eqnarray}
&&\norm{W_\alpha(t)\phi}_{L_x^\infty \widehat L_{t}^r}\les \|\phi\|_{\dot{\widehat{H}}^{-\frac{1+\alpha}{r}}_r}.
\end{eqnarray}
\end{lemma}
\begin{proof}
By the change of variables $\eta=\omega(\xi)$, we get
\begin{equation*}
\begin{aligned}
W_\alpha(t)\phi&=\int e^{it\omega(\xi)}e^{ix\xi}\widehat{\phi}(\xi)d\xi\\
&=\int e^{it\eta}e^{ix\omega^{-1}(\eta)}\widehat{\phi}(\omega^{-1}(\eta))(\omega^{-1}(\eta))^\prime d\eta.
\end{aligned}
\end{equation*}
A straightforward computation with $|\omega^\prime(\xi)|=|\xi|^{1+\alpha}$ gives
\begin{equation*}
\begin{aligned}
\norm{W_\alpha(t)\phi}_{L_x^\infty \widehat L_{t}^r}&=\left\|\widehat{\phi}(\omega^{-1}(\eta))(\omega^{-1}(\eta))^\prime\right\|_{L_{\eta}^{r^\prime}}\\
&=\left(\int \left|\frac{\widehat{\phi}(\xi)}{\omega^\prime(\xi)}\right|^{r^\prime}|\omega^\prime(\xi)|d \xi\right)^{\frac{1}{r^\prime}}\\
&=\left(\int \left||\xi|^{-\frac{1+\alpha}{r}}\widehat{\phi}(\xi)\right|^{r^\prime}d\xi\right)^{\frac{1}{r^\prime}}= \|\phi\|_{\dot{\widehat{H}}^{-\frac{1+\alpha}{r}}_r}.
\end{aligned}
\end{equation*}
\end{proof}

The Duhamel integral of \eqref{eq:dgBO} is given by
\begin{equation}
u(t)=\Gamma u(t):= W_{\alpha}(t)u_0+\int_0^t W_{\alpha}(t-s)\mathcal{N}(u)(s)ds
\end{equation}
where $\mathcal{N}$ is the nonlinear function of $u$. Let $0<T\leq 1$ and $\psi_{T}(t)=\psi(t/T)$. To apply the contraction principle in $\widehat{X}^{s, b}_{r}$, we shall introduce the restriction norm space $\widehat{X}^{s, b}_{r}(T)$ by
\begin{equation}
\|f\|_{\widehat{X}^{s, b}_{r}(T)}:=\inf \left\{\|\tilde{f}\|_{\widehat{X}^{s, b}_{r}}:\left.\tilde{f}\right|_{[-T, T] \times \mathbb{R}}=f\right\}.
\end{equation}
For $u \in \widehat{X}^{s, b}_{r}(T)$ with extension $\tilde{u} \in \widehat{X}^{s, b}_{r}$, an extension of $\Gamma u$ is given by
\begin{equation}
\widetilde{\Gamma u}(t)= \psi(t) W_{\alpha}(t)u_0+\psi_{T}(t)\int_0^t W_{\alpha}(t-s)\mathcal{N}(\tilde{u})(s)ds.
\end{equation}
Thus, it reduces to prove
\begin{equation}
\|\Gamma u\|_{\widehat{X}^{s, b}_{r}(T)} \leq \left\|\psi(t) W_{\alpha}(t)u_0\right\|_{\widehat{X}^{s, b}_{r}}+\left\|\psi_{T}(t)\int_0^t W_{\alpha}(t-s)\mathcal{N}(\tilde{u})(s)ds\right\|_{\widehat{X}^{s, b}_{r}}.
\end{equation}
Furthermore, we prove the following linear estimates.
\begin{proposition}[Linear estimates]\label{proplineares}
(a) Assume $s, b \in \mathbb{R}$, $1 \leq r \leq \infty$ and $\phi \in \widehat{H}^{s}_{r}$. Then there exists $C>0$
such that
\begin{eqnarray}
\norm{\psi(t)W_{\alpha}(t)\phi}_{\widehat{X}^{s, b}_{r}}\leq C\norm{\phi}_{\widehat{H}^{s}_{r}}.
\end{eqnarray}

(b) Assume $s\in \R$, $1<r<\infty$, and $b^\prime+1\geq b\geq 0 \geq b^\prime >-\frac{1}{r^\prime}$. Then there exists $C>0$ such
that
\begin{eqnarray}
\normo{\psi_{T}(t)\int_0^t W_{\alpha}(t-s)F(s)ds}_{\widehat{X}^{s, b}_{r}}\leq
CT^{1+b^\prime-b}\norm{F}_{\widehat{X}^{s, b^\prime}_{r}}.
\end{eqnarray}
\end{proposition}
Therefore, to prove the local well-posedness of  dgBO \eqref{eq:dgBO}, we focus on showing nonlinear estimates in $\widehat{X}^{s, b}_{r}$ spaces. Given $b>\frac{1}{r}$ and $b^\prime \in (b-1,0]$, it turns to prove
\begin{equation}\label{nonlinear}
\|\mathcal{N}(\tilde{u})\|_{\widehat{X}^{s, b^\prime}_{r}} \lesssim \|\tilde{u}\|_{\widehat{X}^{s, b}_{r}}^{2}.
\end{equation}
Without loss of generality, let $b=\frac{1}{r}+\varepsilon$ and $b^{\prime}=-\frac{1}{r^\prime}+2\varepsilon$. For brevity, in the nonlinear estimate \eqref{nonlinear}, we still use $u$ to denote $\tilde{u}$. 
By duality argument, \eqref{nonlinear} is implied as 
\begin{equation}\label{dual_1}
\begin{aligned}
&\left|\int_{\mathbb{R}^{4}} \widehat{u}_{1}(\xi_{1}, \tau_{1})\widehat{u}_{2}(\xi_{2}, \tau_{2})\widehat{v}(\xi_{1}+\xi_{2}, \tau_{1}+\tau_{2})d\xi_{1}d\xi_{2}d\tau_{1}d\tau_{2}\right|\\
\lesssim &\prod^{2}_{j=1}\|u_j\|_{\X_{r}^{s, 1/r+\varepsilon}}\|v\|_{\X_{r^{\prime}}^{-s-1,  1/r^{\prime}-2\varepsilon}}.
\end{aligned}
\end{equation}
For $\xi_1,\xi_2 \in \R$ and $\omega:\R\rightarrow \R$, define the resonance function
\begin{equation}\label{eq:reso}
\Omega(\xi_1,\xi_2)=\omega(\xi_1)+\omega(\xi_2)-\omega(\xi_1+\xi_2),
\end{equation}
which plays a crucial role in bilinear estimates of the $X^{s,b}$-type space. See \cite{T2001} for a perspective discussion. Let $\sigma_j$ and $\sigma$ denote the modulations given by
$$\sigma_i=\tau_i-\omega(\xi_i), \quad  \sigma=\tau-\omega(\xi).$$
Under the restrictions $\xi=\xi_1+\xi_2$ and $\tau=\tau_1+\tau_2$, we have
\begin{equation}\label{L_max}
\tau-\omega(\xi)=\tau_1-\omega(\xi_1)+\tau_2-\omega(\xi_2)+\Omega(\xi_1, \xi_2).
\end{equation}
We apply Littlewood-Paley dyadic decomposition \eqref{dyadic} to each component of \eqref{dual_1}. 
Then to prove \eqref{dual_1}, it suffices to prove
 \begin{equation}\label{dual_2}
\begin{aligned}
&\left|\int_{\mathbb{R}^{4}} \widehat{u}_{1}(\xi_{1}, \tau_{1})\widehat{u}_{2}(\xi_{2}, \tau_{2})\widehat{v}(\xi_{1}+\xi_{2}, \tau_{1}+\tau_{2})d\xi_{1}d\xi_{2}d\tau_{1}d\tau_{2}\right|\\
\lesssim &\sum_{L_1,L_2,L}\sum_{N_1,N_2, N}C(L_1,L_2,L)C(N_1,N_2, N) \\
&\cdot \prod^2_{j=1}\|Q_{L_j} P_{N_j}u_j\|_{\X_{r}^{s, 1/r+\varepsilon}}\|Q_{L} P_{N}v\|_{\X_{r^{\prime}}^{-s-1, 1/r^{\prime}-2\varepsilon}},
\end{aligned}
\end{equation}
where $C(L_1, L_2, L)$ and $C(N_1, N_2, N)$ are suitable bounds that allow us to sum over all dyadic numbers. 

\section{Bilinear Estimates}
In this section, we prove some dyadic bilinear estimates which are crucial for proving \eqref{dual_2} in the next section. 
For compactly supported non-negative functions $f_1, f_2\in L^{r^{\prime}}(\R\times \R)$ with $$\supp(f_i) \subset D_{N_i, L_i}:=\left\{(\xi_i, \sigma_i): |\xi_i|\sim N_i, |\sigma_i|\sim L_i\right\}, \ i=1, 2$$ 
and $f\in L^r(\R\times \R)$ with $$\supp(f) \subset D_{N, L}:=\left\{(\xi, \sigma): |\xi|\sim N, |\sigma|\sim L\right\},$$ 
we define
\begin{equation*}
J\left(f_{1}, f_{2}, f\right)=\int_{\mathbb{R}^{4}} f_{1}\left(\xi_{1}, \sigma_{1}\right) f_{2}\left(\xi_{2}, \sigma_{2}\right) f\left(\xi_{1}+\xi_{2}, \sigma_{1}+\sigma_{2}+\Omega\left(\xi_{1}, \xi_{2}\right)\right) d \xi_{1} d \xi_{2} d \sigma_{1} d \sigma_{2}.
\end{equation*}
It is convenient to define $L_{\max}$, $L_{\med}$ and $L_{\min}$ as the maximum, median and minimum values of $L$, $L_1$ and $L_2$, respectively. Similarly, define $N_{\max}$, $N_{\med}$ and $N_{\min}$ as the maximum, median and minimum values of $N$, $N_1$ and $N_2$.

Now, we state our bilinear estimates.
\begin{lemma}\label{lemma_1}
For any $N_{1}, N_{2}, N$ and $L_{1}, L_{2}, L$, we have
\begin{equation}
J\left(f_{1}, f_{2}, f\right) \lesssim A\cdot B \cdot \left\|f_{1}\right\|_{L^{r^{\prime}}}\left\|f_{2}\right\|_{L^{r^{\prime}}}\left\|f\right\|_{L^{r}}
\end{equation}
where
\begin{equation}
A =L_{\min}^{\frac{1}{r^{\prime}}}\min\left\{L_1^{\frac{1}{r}-\frac{1}{r^{\prime}}}, L_2^{\frac{1}{r}-\frac{1}{r^{\prime}}}\right\},
\end{equation}
and
\begin{equation}
B=N_{\min}^{\frac{1}{r^{\prime}}}\min\left\{N_1^{\frac{1}{r}-\frac{1}{r^{\prime}}}, N_2^{\frac{1}{r}-\frac{1}{r^{\prime}}}\right\}.
\end{equation}
\end{lemma}
\begin{proof}
The proof only relies on the change of variables. We first integrate over $\sigma_{1}$ and $\sigma_{2}$. If $L_2=L_{\max}$, we may use the change of variables $\theta_{1}=\sigma_{1}$ and $\theta_{2}=\sigma_{1}+\sigma_{2}+\Omega\left(\xi_{1}, \xi_{2}\right)$. By H\"older inequality, we get 
\begin{equation*}
\begin{aligned}
& J\left(f_{1}, f_{2}, f\right) \\
=&\int_{\mathbb{R}^{4}} f_{1}\left(\xi_1, \theta_{1}\right) f_{2}\left(\xi_2, \theta_{2}-\theta_{1}-\Omega\left(\xi_1, \xi_2\right) \right)f\left(\xi_1+\xi_2, \theta_{2}\right) d \theta_{1} d \theta_{2} d \xi_1 d \xi_2 \\	\lesssim&\int_{\mathbb{R}^{2}}\left\|f_{1}\left(\xi_1, \theta_{1}\right) f\left(\xi_1+\xi_2, \theta_{2}\right)\right\|_{L_{ \theta_{1}, \theta_{2}}^{r}}\left\|f_{2}\left(\xi_2, \theta_{2}-\theta_{1}-\Omega\left(\xi_1, \xi_2\right) \right)\right\|_{L_{\theta_{1}, \theta_{2}}^{r^{\prime}}} d \xi_1 d \xi_2\\
\lesssim&L_1^{\frac{1}{r}-\frac{1}{r^{\prime}}}\min\left\{L_1^{\frac{1}{r^{\prime}}}, L^{\frac{1}{r^{\prime}}}\right\}\int_{\mathbb{R}^{2}}\left\|f_{1}\left(\xi_1, \sigma_1\right)\right\|_{L_{\sigma_1}^{r^\prime}}\left\|f\left(\xi_1+\xi_2, \sigma\right)\right\|_{L_{\sigma}^{r}}\left\|f_{2}\left(\xi_2, \sigma_2\right)\right\|_{L_{\sigma_2}^{r^{\prime}}} d \xi_1 d \xi_2.
\end{aligned}
\end{equation*}
Then, the integration over $\xi_1$ and $\xi_2$ 
\begin{equation}\label{case_1}
\int_{\mathbb{R}^{2}}\left\|f_{1}\left(\xi_1, \sigma_1\right)\right\|_{L_{\sigma_1}^{r^\prime}}\left\|f\left(\xi_1+\xi_2, \sigma\right)\right\|_{L_{\sigma}^{r}}\left\|f_{2}\left(\xi_2, \sigma_2\right)\right\|_{L_{\sigma_2}^{r^{\prime}}} d \xi_1 d \xi_2
\end{equation}
pushes us to consider three subcases as follows.
\begin{itemize} 
\item Subcase 1.1. We may integrate $f_1$ and $f_2$ first. By H\"older inequality, we have
\end{itemize}
\begin{equation*}
\eqref{case_1}\lesssim \left\|f_{1}\right\|_{L^{r^\prime}}\left\|f_{2}\right\|_{L^{r^{\prime}}}\min\left\{N_{1}^{\frac{1}{r}}, N_{2}^{\frac{1}{r}}\right\}\left\|f\right\|_{L^{r}}.
\end{equation*}
\begin{itemize} 
\item Subcase 1.2. We also may integrate $f_1$ and $f$ first. A change of variables $\ell=\xi_1$ and $\eta=\xi_1+\xi_2$ gives
\end{itemize}
\begin{equation*}
\begin{aligned}
\eqref{case_1}=& \int_{\mathbb{R}^{2}}\left\|f_{1}\left(\ell, \sigma_1\right)\right\|_{L_{\sigma_1}^{r^\prime}}\left\|f\left(\eta, \sigma\right)\right\|_{L_{\sigma}^{r}}\left\|f_{2}\left(\eta-\ell, \sigma_2\right)\right\|_{L_{\sigma_2}^{r^{\prime}}} d \eta d \ell\\
\lesssim& N_1^{\frac{1}{r}-\frac{1}{r^{\prime}}}\left\|f_{1}\right\|_{L^{r^\prime}}\left\|f\right\|_{L^{r}}\min\left\{N_1^{\frac{1}{r^{\prime}}}, N^{\frac{1}{r^{\prime}}}\right\}\left\|f_{2}\right\|_{L^{r^{\prime}}}.
\end{aligned}
\end{equation*}
\begin{itemize} 
\item Subcase 1.3. Lastly, we may integrate $f_2$ and $f$ first. A change of variables $\ell=\xi_2$ and $\eta=\xi_1+\xi_2$ gives
\end{itemize}
\begin{equation*}
\begin{aligned}
\eqref{case_1}=& \int_{\mathbb{R}^{2}}\left\|f_{1}\left(\eta-\ell, \sigma_1\right)\right\|_{L_{\sigma_1}^{r^\prime}}\left\|f\left(\eta, \sigma\right)\right\|_{L_{\sigma}^{r}}\left\|f_{2}\left(\ell, \sigma_2\right)\right\|_{L_{\sigma_2}^{r^{\prime}}} d \eta d \ell\\
\lesssim& N_2^{\frac{1}{r}-\frac{1}{r^{\prime}}}\left\|f_{2}\right\|_{L^{r^\prime}}\left\|f\right\|_{L^{r}}\min\left\{N_2^{\frac{1}{r^{\prime}}}, N^{\frac{1}{r^{\prime}}}\right\}\left\|f_{1}\right\|_{L^{r^{\prime}}}.
\end{aligned}
\end{equation*}
Thus, if $L_2=L_{\max}$, we may obtain
\begin{equation}
\begin{aligned}
J\left(f_{1}, f_{2}, f\right)\lesssim B \cdot L_1^{\frac{1}{r}-\frac{1}{r^{\prime}}}\min\left\{L_1^{\frac{1}{r^{\prime}}}, L^{\frac{1}{r^{\prime}}}\right\}\cdot \left\|f_{1}\right\|_{L^{r^{\prime}}}\left\|f_{2}\right\|_{L^{r^{\prime}}}\left\|f\right\|_{L^{r}},
\end{aligned}
\end{equation}
where $B=N_{\min}^{\frac{1}{r^{\prime}}}\min\left\{N_1^{\frac{1}{r}-\frac{1}{r^{\prime}}}, N_2^{\frac{1}{r}-\frac{1}{r^{\prime}}}\right\}$.

$f_1$ and $f_2$ are symmetric, hence for $L_1=L_{\max}$ we have 
\begin{equation}
\begin{aligned}
J\left(f_{1}, f_{2}, f\right)\lesssim B \cdot L_2^{\frac{1}{r}-\frac{1}{r^{\prime}}}\min\left\{L_2^{\frac{1}{r^{\prime}}}, L^{\frac{1}{r^{\prime}}}\right\}\cdot \left\|f_{1}\right\|_{L^{r^{\prime}}}\left\|f_{2}\right\|_{L^{r^{\prime}}}\left\|f\right\|_{L^{r}}.
\end{aligned}
\end{equation}
Similarly, it is easy to verify that the case $L=L_{\max}$ implies
\begin{equation}
J\left(f_{1}, f_{2}, f\right) \lesssim B \cdot \min \left\{L_{1}^{\frac{1}{r}}, L_{2}^{\frac{1}{r}}\right\} \cdot \left\|f_{1}\right\|_{L^{r^{\prime}}}\left\|f_{2}\right\|_{L^{r^{\prime}}}\left\|f\right\|_{L^{r}}. 
\end{equation}
This completes the proof.
\end{proof}

In the sequel, we will consider different interactions. To take advantage of the resonance function $\Omega(\xi_1, \xi_2)$, we shall first integrate over $\xi_1$ and $\xi_2$. 
\begin{lemma}\label{HL} ($high \times low$)

(a) If $L=L_{\max}$, we have
\begin{equation}\label{HL_L}
J\left(f_{1}, f_{2}, f\right) \lesssim L_{1}^{\frac{1}{r}}L_{2}^{\frac{1}{r}}N_{\max}^{-\frac{1+\alpha}{r}}\left\|f_{1}\right\|_{L^{r^{\prime}}}\left\|f_{2}\right\|_{L^{r^{\prime}}}\left\|f\right\|_{L^{r}}.
\end{equation}

(b) If $L_i=L_{\max}\ (i=1, 2)$, we have
\begin{equation}\label{HL_L12}
J\left(f_{1}, f_{2}, f\right) \lesssim L_{1}^{\frac{1}{r}}L_{2}^{\frac{1}{r}}L^{\frac{1}{r^\prime}}L_{i}^{-\frac{1}{r}}N_{\min}^{\frac{1}{r}-\frac{1}{r^{\prime}}}(N_{\max}^{\alpha}N_{i})^{-\frac{1}{r^\prime}}\left\|f_{1}\right\|_{L^{r^{\prime}}}\left\|f_{2}\right\|_{L^{r^{\prime}}}\left\|f\right\|_{L^{r}}.
\end{equation}
\end{lemma}

\begin{proof}
Without loss of generality, we assume $N_1 \gg N_2$. Note that $\Omega(\xi_{1}, \xi_{2})=\omega(\xi_{1})+\omega(\xi_{2})-\omega(\xi_{1}+\xi_{2})$, where $\omega(\xi)=-\xi|\xi|^{1+\alpha}$ and $|\omega^{\prime}(\xi)|\sim|\xi|^{1+\alpha}$.

If $L=L_{\max}$, we may integrate $f_1$ and $f_2$ together, which gives
\begin{equation}\label{HL_1}
J\left(f_{1}, f_{2}, f\right)
\lesssim L_{1}^{\frac{1}{r}}L_{2}^{\frac{1}{r}}\left\|f_{1} \right\|_{L^{r^{\prime}}}\left\|f_{2}\right\|_{L^{r^{\prime}}}\left\|f\left(\xi_{1}+\xi_{2}, \Omega\left(\xi_{1}, \xi_{2}\right)\right)\right\|_{L_{ \xi_{1}, \xi_{2}}^{r}}.
\end{equation}
Let $\eta=\xi_{1}+\xi_{2}$ and $\ell=\Omega\left(\xi_{1}, \xi_{2}\right)$. The Jacobian is 
\begin{equation*}
\frac{\partial(\eta, \ell)}{\partial\left(\xi_{1}, \xi_{2}\right)}=
\left|\begin{array}{cc}
1 & 1 \\
\partial_{\xi_{1}} \Omega & \partial_{\xi_{2}} \Omega
\end{array}\right|
=\partial_{\xi_{2}} \Omega-\partial_{\xi_{1}} \Omega =\omega^{\prime}(\xi_{2})-\omega^{\prime}(\xi_{1}).
\end{equation*}
In addition, in $high \times low$ interactions, $|\xi_1| \gg |\xi_2|$ yields 
\begin{equation}
\left|\omega^{\prime}(\xi_{1})-\omega^{\prime}(\xi_{2})\right|\sim \left||\xi_1|^{1+\alpha}-|\xi_2|^{1+\alpha}\right|\sim|\xi_1|^{1+\alpha}.
\end{equation}
Therefore, 
\begin{equation*}
\begin{aligned}
\eqref{HL_1}&\lesssim L_{1}^{\frac{1}{r}}L_{2}^{\frac{1}{r}}\left\|f_{1}\right\|_{L^{r^{\prime}}}\left\|f_{2}\right\|_{L^{r^{\prime}}}\left(\int_{\mathbb{R}^{2}}\left|f(\eta, \ell)\right|^{r} \frac{1}{\mid \omega^{\prime}\left(\xi_{1}\right)-\omega^{\prime} \left(\xi_{2}\right) \mid} d \eta d \ell\right)^{\frac{1}{r}}\\
&\lesssim  L_{1}^{\frac{1}{r}}L_{2}^{\frac{1}{r}}\left\|f_{1} \right\|_{L^{r^{\prime}}}\left\|f_{2}\right\|_{L^{r^{\prime}}}(N_{1}^{1+\alpha})^{-\frac{1}{r}}\left\|f\right\|_{L^{r}}.
\end{aligned}
\end{equation*}
We complete part (a).

If $L_{2}=L_{\max}$, we may integrate $f_1$ and $f$ together. The change of variables $\theta_1=\sigma_1$ and $\theta_2=\sigma_1+\sigma_2+\Omega\left(\xi_1, \xi_2\right)$ gives
\begin{equation}\label{HL_2}
\begin{aligned}
&J\left(f_{1}, f_{2}, f\right) \\
=& \int_{\mathbb{R}^{4}} f_{1}\left(\xi_{1}, \theta_{1}\right) f_{2}\left(\xi_{2}, \theta_2-\theta_1-\Omega\left(\xi_1, \xi_2\right)\right) f\left(\xi_{1}+\xi_{2}, \theta_2\right) d \xi_{1} d \xi_{2} d \theta_{1} d \theta_{2}\\
\lesssim & L_{1}^{\frac{1}{r}}L^{\frac{1}{r^{\prime}}}N_{2}^{\frac{1}{r}-\frac{1}{r^{\prime}}}\left\|f_{1}\right\|_{L^{r^{\prime}}}\left\|f\right\|_{L^{r}} \left\|f_{2}\left(\xi_{2}, -\Omega\left(\xi_1, \xi_2\right)\right)\right\|_{L_{\xi_1, \xi_{2}}^{r^{\prime}}}\\
\lesssim &L_{1}^{\frac{1}{r}}L^{\frac{1}{r^{\prime}}}N_{2}^{\frac{1}{r}-\frac{1}{r^{\prime}}}\left\|f_{1}\right\|_{L^{r^{\prime}}}\left\|f\right\|_{L^{r}} \left(\int_{\mathbb{R}^{2}}\left|f_{2}(\xi_2, \phi)\right|^{r^{\prime}} \left|\frac{d \xi_1}{d \phi} \right| d \phi d \xi_2\right)^{\frac{1}{r^{\prime}}}\\
\lesssim & L_{1}^{\frac{1}{r}}L^{\frac{1}{r^{\prime}}}N_{2}^{\frac{1}{r}-\frac{1}{r^{\prime}}}\left(N_{1}^{\alpha}N_{2}\right)^{-\frac{1}{r^{\prime}}}\left\|f_{1}\right\|_{L^{r^{\prime}}}\left\|f\right\|_{L^{r}}\left\|f_{2}\right\|_{L^{r^{\prime}}}.
\end{aligned}
\end{equation}
Here, the Jacobian is 
\begin{equation*}
    \begin{aligned}
        \left|\frac{d \phi}{d \xi_1}\right|=\left|\partial_{\xi_1}\Omega(\xi_1, \xi_2)\right|
        =\left|\omega^{\prime}(\xi_{1})-\omega^{\prime}(\xi_1+\xi_{2})\right|. 
    \end{aligned}
\end{equation*}
Note that $\left|\xi_{1}\right| \gg\left|\xi_{2}\right|$. It means that $\xi_1$ and $\xi_1+\xi_2$ have the same signs. By the mean value theorem, we have
\begin{equation}
\left|\omega^{\prime}(\xi_{1})-\omega^{\prime}(\xi_1+\xi_{2})\right| \sim \left|\xi_1\right|^{\alpha}\left|\xi_2\right|\sim N_{1}^{\alpha}N_{2}.
\end{equation}

If $L_{1}=L_{\max}$, similarly, we have
\begin{equation}
\begin{aligned}
&J\left(f_{1}, f_{2}, f\right) \\ =& \int_{\mathbb{R}^{4}} f_{1}\left(\xi_{1}, \theta_1-\theta_2-\Omega\left(\xi_1, \xi_2\right)\right) f_{2}\left(\xi_{2}, \theta_{2}\right) f\left(\xi_{1}+\xi_{2}, \theta_1\right) d \xi_{1} d \xi_{2} d \theta_{1} d \theta_{2}\\
\lesssim & L_{2}^{\frac{1}{r}}L^{\frac{1}{r^{\prime}}}N_{2}^{\frac{1}{r}-\frac{1}{r^{\prime}}}\left\|f_{2}\right\|_{L^{r^{\prime}}}\left\|f\right\|_{L^{r}} \left(\int_{\mathbb{R}^{2}}\left|f_{1}(\xi_1, \phi)\right|^{r^{\prime}} \left|\frac{d \xi_2}{d \phi} \right| d \phi d \xi_1\right)^{\frac{1}{r^{\prime}}}\\
\lesssim & L_{2}^{\frac{1}{r}}L^{\frac{1}{r^{\prime}}}N_{2}^{\frac{1}{r}-\frac{1}{r^{\prime}}}\left\|f_{2}\right\|_{L^{r^{\prime}}}\left\|f\right\|_{L^{r}} \left(N_{1}^{1+\alpha}\right)^{-\frac{1}{r^{\prime}}}\left\|f_{1}\right\|_{L^{r^{\prime}}},
\end{aligned}
\end{equation}
where we use the fact $|\xi_1+\xi_2|\gg |\xi_2|$ and obtain
\begin{equation}
\begin{aligned}
\left|\frac{d \phi}{d \xi_2}\right|=\left|\partial_{\xi_2}\Omega(\xi_1, \xi_2)\right|
=\left|\omega^{\prime}(\xi_{2})-\omega^{\prime}(\xi_1+\xi_{2})\right|\sim \left|\xi_1\right|^{1+\alpha}\sim N_{1}^{1+\alpha}. 
\end{aligned}
\end{equation}
Part (b) is completed.
\end{proof}

In what follows, we consider $high \times high$ interactions, which indicates $N_1 \sim N_2$. In particular, when the frequencies of the two factors are very close, ($\left|\xi_{1}-\xi_{2}\right| \ll |\xi_1| \sim |\xi_2|$), we go back to the proof of Lemma 9 in \cite{Gr2010}.

\begin{lemma} \label{HH} ($high \times high$)

(a) Assume $L=L_{\max}$. If $\xi_1\xi_2<0$, we have
\begin{equation}\label{HHL}
J\left(f_{1}, f_{2}, f\right) \lesssim L_{1}^{\frac{1}{r}}L_{2}^{\frac{1}{r}}(N_{\max}^{\alpha}N)^{-\frac{1}{r}}\left\|f_{1}\right\|_{L^{r^{\prime}}}\left\|f_{2}\right\|_{L^{r^{\prime}}}\left\|f\right\|_{L^{r}}.
\end{equation}
If $\xi_1\xi_2>0$, we have for $b>\frac{1}{r}$
\begin{equation}\label{HHH}
\begin{aligned}
J\left(f_{1}, f_{2}, f\right)
\lesssim N_{\max}^{-\frac{\alpha}{2r}}\left\|\mathcal{F}^{-1}f_{1}\right\|_{\X_r^{0, b}}\left\|\mathcal{F}^{-1}f_{2}\right\|_{\X_r^{0, b}}\left\|f\right\|_{L^{r}}.
\end{aligned}
\end{equation}

(b) If $L_i=L_{\max} \ (i=1,2)$, we have
\begin{equation}\label{HH_L12}
J\left(f_{1}, f_{2}, f\right) \lesssim L_{1}^{\frac{1}{r}}L_{2}^{\frac{1}{r}}L^{\frac{1}{r^\prime}}L_{j}^{-\frac{1}{r}}N_{\max}^{\frac{1}{r}-\frac{1}{r^{\prime}}}(N_{\max}^{1+\alpha})^{-\frac{1}{r^\prime}}\left\|f_{1}\right\|_{L^{r^{\prime}}}\left\|f_{2}\right\|_{L^{r^{\prime}}}\left\|f\right\|_{L^{r}}.
\end{equation}
\end{lemma}

\begin{proof}
If $\xi_1$ and $\xi_2$ have different signs ($\xi_1\xi_2<0$), we get $N_1 \sim N_2 \gg N$ or $N_1 \sim N_2 \sim N$. It implies that $\xi_1$ and $\xi_1-\xi$ have the same signs. By the mean value theorem, we have
\begin{equation}
\left|\omega^{\prime}(\xi_{1})-\omega^{\prime}(\xi_{2})\right|=\left|\omega^{\prime}(\xi_{1})-\omega^{\prime}(\xi-\xi_1)\right|\sim \left|\xi_1\right|^{\alpha}\left|\xi\right|\sim N_{1}^{\alpha}N.
\end{equation}
Hence, with $L=L_{\max}$, we get
\begin{equation}\label{HHL_2}
\begin{aligned}
J\left(f_{1}, f_{2}, f\right) 
\lesssim &L_{1}^{\frac{1}{r}}L_{2}^{\frac{1}{r}}\left\|f_{1} \right\|_{L^{r^{\prime}}}\left\|f_{2}\right\|_{L^{r^{\prime}}}\left\|f\right\|_{L^{r}} \mid \omega^{\prime}(\xi_{1})-\omega^{\prime} (\xi_{2}) \mid^{-\frac{1}{r}}\\
\lesssim & L_{1}^{\frac{1}{r}}L_{2}^{\frac{1}{r}}(N_{1}^{\alpha}N)^{-\frac{1}{r}}\left\|f_{1} \right\|_{L^{r^{\prime}}}\left\|f_{2}\right\|_{L^{r^{\prime}}}\left\|f\right\|_{L^{r}}.
\end{aligned}
\end{equation}
The desired \eqref{HHL} is proved.

If $\xi_1$ and $\xi_2$ have the same signs ($\xi_1\xi_2>0$), we need to perform a delicate discussion on $|\xi_{1}-\xi_2|$. Hence, we decompose dyadically with respect to $|\xi_{1}-\frac{\xi}{2}| =\frac{1}{2} |\xi_1-\xi_2|$ and discuss contributions on $A_k$, where
\begin{equation}
A_k=\{\xi_1: |\xi_{1}-\xi/2|\sim 2^{-k}, k \geq 1 \} \ \text{and} \ A_0=\{\xi_1: 1 \leq |\xi_{1}-\xi/2|\ll |\xi_1|\}.
\end{equation}
In this subcase, we go back to the following estimate
\begin{equation*}
\begin{aligned}
J\left(f_{1}, f_{2}, f\right)
\lesssim  \left\|\int_{*} f_{1}\left(\xi_{1}, \tau_{1}\right)f_{2}\left(\xi_2, \tau_2\right)d\xi_1 d\tau_1\right\|_{L_{\xi, \tau}^{r^{\prime}}}\left\|f\right\|_{L_{\xi,\tau}^{r}},
\end{aligned}
\end{equation*}
and focus on  
\begin{equation}\label{A_k}
\left\|\int_{*} \chi_{A_k}(\xi_1) f_{1}\left(\xi_{1}, \tau_{1}\right)f_{2}\left(\xi_2, \tau_2\right)d\xi_1 d\tau_1\right\|_{L_{\xi, \tau}^{r^{\prime}}},
\end{equation}
where $\int_{*}$ is shorthand for $\int_{\xi=\xi_1+\xi_2, \tau=\tau_1+\tau_2}$. For $i=1, 2$, we choose $g_i$ with $\|g_i\|_{L_{\xi, \tau}^{r^{\prime}}}=\|\mathcal{F}^{-1}f_i\|_{\X_r^{0, b}}$ so that \eqref{A_k} is rewritten as
\begin{equation*}\label{g_i}
\left\|\int_{*} \chi_{A_k}(\xi_1) \frac{g_{1}\left(\xi_{1}, \tau_{1}\right)g_{2}\left(\xi_2, \tau_2\right)}{\langle \sigma_1 \rangle ^{b}\langle \sigma_2 \rangle ^{b}}d\xi_1 d\tau_1\right\|_{L_{\xi, \tau}^{r^{\prime}}}.
\end{equation*}
We use Hölder's inequality over $\tau_1$ and Lemma 4.2 in \cite{GTV1997} to obtain
\begin{equation*}
\int_{*} \frac{g_{1}\left(\xi_{1}, \tau_{1}\right)g_{2}\left(\xi_2, \tau_2\right)}{\langle \sigma_1 \rangle ^{b}\langle \sigma_2 \rangle ^{b}} d\tau_1 
\lesssim \langle\sigma_{r e s}\rangle^{-b}\left(\int_{*} \left|g_1\left(\xi_1, \tau_1\right) g_2\left(\xi_2, \tau_2\right)\right|^{r^{\prime}}d \tau_1\right)^{\frac{1}{r^{\prime}}},
\end{equation*}
where $$\sigma_{res}=\tau-\xi_1|\xi_1|^{1+\alpha}-\xi_2|\xi_2|^{1+\alpha}=\tau-2^{-1-\alpha}\xi|\xi|^{1+\alpha}-h(x)$$ 
with $x=\xi_1-\xi/2$ and 
$$h(x)=(\xi/2+x)|\xi/2+x|^{1+\alpha}+(\xi/2-x)|\xi/2-x|^{1+\alpha}-2^{-1-\alpha}\xi|\xi|^{1+\alpha}.$$
Note that $h(0)=0$. Moreover,
$$h^{\prime}(x) \sim x|\xi|^{\alpha} \ \text{and} \ h^{\prime}(0)=0.$$ 

Now we turn to the case $k=0$. Applying Hölder's inequality over $\xi_1$ gives
\begin{equation*}
\begin{aligned}
&\int_{*} \chi_{A_0}(\xi_1) |\xi_1-\xi/2|^{\frac{1}{r}}\frac{g_{1}\left(\xi_{1}, \tau_{1}\right)g_{2}\left(\xi_2, \tau_2\right)}{\langle \sigma_1 \rangle ^{b}\langle \sigma_2 \rangle ^{b}}|\xi_1-\xi/2|^{-\frac{1}{r}}d\xi_1 d\tau_1\\
\lesssim & \left(\int_{A_0}|\xi_1-\xi/2|\langle\sigma_{r e s}\rangle^{-br} d\xi_1\right)^{\frac{1}{r}}\left(\int_{*} |\xi_1-\xi/2|^{-\frac{r^\prime}{r}}\left|g_1\left(\xi_1, \tau_1\right) g_2\left(\xi_2, \tau_2\right)\right|^{r^{\prime}}d\xi_1d \tau_1\right)^{\frac{1}{r^{\prime}}}.
\end{aligned}
\end{equation*}
For the first factor, by the change of variables we have
$$\left(\int \langle\tau-2^{-1-\alpha}\xi|\xi|^{1+\alpha}-h\rangle^{-br} |\xi|^{-\alpha} dh\right)^{\frac{1}{r}}\lesssim |\xi|^{-\frac{\alpha}{r}}.$$ 
For the second one, using Fubini we arrive at
\begin{equation*}
\begin{aligned}
&\left\| \left(\int_{*} |\xi_1-\xi/2|^{-\frac{r^\prime}{r}}\left|g_1\left(\xi_1, \tau_1\right) g_2\left(\xi_2, \tau_2\right)\right|^{r^{\prime}}d\xi_1d \tau_1\right)^{\frac{1}{r^{\prime}}}\right\|_{L_\tau^{r^{\prime}}} \\
\lesssim & \left(\int_* |\xi_1-\xi/2|^{-\frac{r^\prime}{r}}\left\|g_1\left(\xi_1, \cdot\right)\right\|_{L_\tau^{r^{\prime}}}^{r^{\prime}}\left\|g_2\left(\xi_2, \cdot\right)\right\|_{L_\tau^{r^{\prime}}}^{r^{\prime}}d\xi_1\right)^{\frac{1}{r^{\prime}}} \\
\lesssim & \left(\int_* \left\|g_1\left(\xi_1, \cdot\right)\right\|_{L_\tau^{r^{\prime}}}^{r^{\prime}}\left\|g_2\left(\xi_2, \cdot\right)\right\|_{L_\tau^{r^{\prime}}}^{r^{\prime}}d\xi_1\right)^{\frac{1}{r^{\prime}}}.
\end{aligned}
\end{equation*}
Taking $\|\ \|_{L_{\xi}^{r^{\prime}}}$ gives the bound for the contribution from $A_0$.

Then we turn to the case $k\geq 1$. By Taylor's expansion with $|x|<1$, we obtain $|h(x)| \sim |x|^{2}|\xi|^{\alpha}$ and another estimate of $|h^{\prime}(x)|$; that is, $ |h^{\prime}(x) |\sim |x||\xi|^{\alpha} \sim |h(x)|^{\frac{1}{2}}|\xi|^{\frac{\alpha}{2}}$. Therefore, we have
\begin{equation*}
\begin{aligned}
&\int_{*} \chi_{A_k}(\xi_1) \frac{g_{1}\left(\xi_{1}, \tau_{1}\right)g_{2}\left(\xi_2, \tau_2\right)}{\langle \sigma_1 \rangle ^{b}\langle \sigma_2 \rangle ^{b}}d\xi_1 d\tau_1\\
\lesssim & \left(\int_{  \chi_{A_k}}\langle\sigma_{r e s}\rangle^{-b r}d \xi_1\right)^{\frac{1}{r}}\left(\int_* \left|g_1\left(\xi_1, \tau_1\right) g_2\left(\xi_2, \tau_2\right)\right|^{r^{\prime}}d \xi_1 d \tau_1 \right)^{\frac{1}{r^{\prime}}}.
\end{aligned}
\end{equation*}
For the first factor, we now use $dh=|h|^{\frac{1}{2}}|\xi|^{\frac{\alpha}{2}}dx$ and obtain
\begin{equation}
\left(\int_{\chi_{A_k}}\langle\sigma_{r e s}\rangle^{-b r}d \xi_1\right)^{\frac{1}{r}} \lesssim \left( \int \langle \tau-2^{-1-\alpha}\xi|\xi|^{1+\alpha}-h\rangle^{-b r}|h|^{-\frac{1}{2}}|\xi|^{-\frac{\alpha}{2}}dh\right)^{\frac{1}{r}} \lesssim |\xi|^{-\frac{\alpha}{2r}}.
\end{equation}
Taking $\|\ \|_{L_{\xi, \tau}^{r^{\prime}}}$ for the second factor gives the bound for $k\geq 1$. Since $|\xi|\sim N_{\max}$, then the desired \eqref{HHH} is proved. Part (a) is completed.

If $L_i=L_{\max} \ (i=1, 2)$, the proof follows a similar argument as \eqref{HL_L12} in Lemma \ref{HL}, where the Jacobian is $N_1^{1+\alpha}$. Hence, we complete part (b).
\end{proof}

\section{Proof of Theorem \ref{theorem}}
In this section, we will prove \eqref{dual_2} case by case, which is sufficient for the proof of Theorem \ref{theorem}. The main ingredients are dyadic bilinear estimates in Section 3 and the use of resonance function \eqref{eq:reso}. We still assume $N_1\geq N_2$. In addition, \eqref{L_max} implies 
\begin{equation}\label{Lmax}
L_{\max}\gtrsim |\Omega\left(\xi_{1}, \xi_{2}\right)|\gtrsim N_{\max}^{\alpha}N_{\min}.
\end{equation}
To simplify notations, we still write $\widehat{u}_{i}:=\mathcal{F}(Q_{L_i} P_{N_i}u_i)(\xi_{i}, \tau_{i})$ for $i=1,2$, and $\widehat{v}:=\mathcal{F}(Q_{L}P_Nv)(\xi, \tau)$. Moreover, we are able to check the boundedness of $C(L_1, L_2, L)$ is easily satisfied. Hence, it suffices to discuss the following boundedness
\begin{equation}
\left|J(\hat{u}_{1}, \hat{u}_{2}, \hat{v})\right|\lesssim C(N_1,N_2, N)\prod^2_{i=1}\|u_i\|_{\X_{r}^{s, 1/r+\varepsilon}}\|v\|_{\X_{r^{\prime}}^{-s-1, 1/r^{\prime}-2\varepsilon}}.
\end{equation}

\begin{proof}[Proof of Theorem \ref{theorem}]
{\bf Case 1:} $N_1\lesssim 1$  or  $ N_1 \gg 1 $ and $L_{\max} \gg N_1^{9}$. These two cases are straightforward and can be solved easily by using Lemma \ref{lemma_1}. Therefore, we can conclude the following.

{$\bullet$} $N_1\lesssim 1$. Since $|\xi_{1}|\lesssim 1$, then $\langle \xi_1 \rangle \sim 1$. 
\begin{align*}
\left|J(\hat{u}_{1}, \hat{u}_{2}, \hat{v})\right|
\lesssim& A \cdot B \cdot \left(L_1L_2\right)^{-\frac{1}{r}-}L^{-\frac{1}{r^\prime}++}\|u_{1}\|_{\X_r^{s, 1/r+}}\|u_{2}\|_{\X_r^{s, 1/r+}}\|v \|_{\X_{r^{\prime}}^{-s-1, 1/r^{\prime}--}},
\end{align*}
which is summable, where $A, B$ are defined in Lemma \ref{lemma_1}.
	
{$\bullet$}$ N_1 \gg 1 $ and $L_{\max} \gg N_1^{9}$.
We may assume $L_1=L_{\max}$, and the other cases follow by a similar argument. We see that
\begin{align*}
\left|J(\hat{u}_{1}, \hat{u}_{2}, \hat{v})\right|
\lesssim & L_{2}^{\frac{1}{r}-\frac{1}{r^\prime}}\min \left\{L_{2}^{\frac{1}{r^\prime}}, L^{\frac{1}{r^\prime}}\right\} N_{1}^{\frac{1}{r}}\|\hat{u}_{1}\|_{L^{r^{\prime}}} \|\hat{u}_{2}\|_{L^{r^{\prime}}}  \|\hat{v}\|_{L^{r}} \\
\lesssim & L_{2}^{\frac{1}{r}-\frac{1}{r^\prime}}\min \left\{L_{2}^{\frac{1}{r^\prime}}, L^{\frac{1}{r^\prime}}\right\} N_{1}^{\frac{1}{r}}N_{1}^{-\frac{9}{r}-}N_{1}^{-s}L_{2}^{-\frac{1}{r}-}N_{2}^{-s}L^{-\frac{1}{r^{\prime}}++}N^{s+1}\\
&\cdot\|u_{1}\|_{\X_r^{s, 1/r+}}\|u_{2}\|_{\X_r^{s, 1/r+}}\|v \|_{\X_{r^{\prime}}^{-s-1, 1/r^{\prime}--}}\\
\lesssim& N_{1}^{-\theta}\prod^2_{i=1}\|u_{i}\|_{\X_{r}^{s, 1/r+}}\|v\|_{\X_{r^{\prime}}^{-s-1, 1/r^{\prime}--}}\\
\end{align*}
which is bounded for some $\theta>0$. 
	
Therefore, we assume that  
\begin{equation}
N_1 \gg 1 \quad \text{and} \quad L_{\max} \les N_1^{9}
\end{equation}
in the following.

{\bf Case 2:} $N \sim N_{1}\gg N_{2}$. In this case, the resonant function $|\Omega\left(\xi_{1}, \xi_{2}\right)|\gtrsim |\xi_1|^{1+\alpha}|\xi-\xi_{1}|\sim N_1^{1+\alpha}N_{2}$.

{$\bullet$} Subcase 2a: $N_2 \lesssim 1$. Applying \eqref{HL_L} directly, we have
\begin{align*}
\left|J(\hat{u}_{1}, \hat{u}_{2}, \hat{v})\right|
\lesssim & L_{1}^{\frac{1}{r}}L_{2}^{\frac{1}{r}}N_{1}^{-\frac{1+\alpha}{r}}\|\hat{u}_{1}\|_{L^{r^{\prime}}} \|\hat{u}_{2}\|_{L^{r^{\prime}}}  \|\hat{v}\|_{L^{r}} \\
\lesssim &L_{1}^{\frac{1}{r}}L_{2}^{\frac{1}{r}}N_{1}^{-\frac{1+\alpha}{r}}L_{1}^{-\frac{1}{r}-}N_{1}^{-s}L_{2}^{-\frac{1}{r}-}L^{-\frac{1}{r^{\prime}}++}N^{s+1}\\
&\cdot\|u_{1}\|_{\X_r^{s, 1/r+}}\|u_{2}\|_{\X_r^{s, 1/r+}}\|v\|_{\X_{r^{\prime}}^{-s-1, 1/r^{\prime}--}}\\
\lesssim& N_{1}^{1-\frac{1+\alpha}{r}}\|u_{1}\|_{\X_r^{s, 1/r+}}\|u_{2}\|_{\X_r^{s, 1/r+}}\|v\|_{\X_{r^{\prime}}^{-s-1, 1/r^{\prime}--}}.
\end{align*}
We find the condition $r<1+\alpha$, which is our first restriction for $r$, depending on $\alpha$. Then in what follows, we assume $N_2 \gg 1$ so that $|\Omega\left(\xi_{1}, \xi_{2}\right)|$ can contribute.

{$\bullet$} Subcase 2b: $L=L_{\max}$. By applying \eqref{HL_L} and \eqref{Lmax}, we get
\begin{align*}
\left|J(\hat{u}_{1}, \hat{u}_{2}, \hat{v})\right|
\lesssim & L_{1}^{\frac{1}{r}}L_{2}^{\frac{1}{r}}N_{1}^{-\frac{1+\alpha}{r}}\|\hat{u}_{1}\|_{L^{r^{\prime}}} \|\hat{u}_{2}\|_{L^{r^{\prime}}}  \|\hat{v}\|_{L^{r}} \\
\lesssim &L_{1}^{\frac{1}{r}}L_{2}^{\frac{1}{r}}N_{1}^{-\frac{1+\alpha}{r}}L_{1}^{-\frac{1}{r}-}N_{1}^{-s}L_{2}^{-\frac{1}{r}-}N_{2}^{-s}L^{-\frac{1}{r^{\prime}}++}N^{s+1}\\
&\cdot\|u_{1}\|_{\X_r^{s, 1/r+}}\|u_{2}\|_{\X_r^{s, 1/r+}}\|v\|_{\X_{r^{\prime}}^{-s-1, 1/r^{\prime}--}}\\
\lesssim & N_{1}^{-\frac{1+\alpha}{r}}N_{1}^{-s}N_{2}^{-s}(N_1^{1+\alpha}N_{2})^{-\frac{1}{r^{\prime}}}N_1^{+}N^{s+1}\\
&\cdot\|u_{1}\|_{\X_r^{s, 1/r+}}\|u_{2}\|_{\X_r^{s, 1/r+}}\|v\|_{\X_{r^{\prime}}^{-s-1, 1/r^{\prime}--}}\\
\lesssim& N_{1}^{-\alpha+}N_{2}^{-s-\frac{1}{r^{\prime}}}\|u_{1}\|_{\X_r^{s, 1/r+}}\|u_{2}\|_{\X_r^{s, 1/r+}}\|v\|_{\X_{r^{\prime}}^{-s-1, 1/r^{\prime}--}}.
\end{align*}
We end up with the conditions $\alpha>0$ and $s > -1-\alpha+\frac{1}{r}$. Note that $\alpha$ here is restricted to $\alpha>0$. Hence, the result fails to contain the endpoint $\alpha=0$.
\begin{remark}
Since $L=L_{\max} \les N_1^{9}$, we have 
\begin{equation}
L^{-\frac{1}{r^{\prime}}++}\lesssim L^{-\frac{1}{r^{\prime}}-}N_1^{+},
\end{equation}
which is used in the third step. We will use the same strategies in the following.
\end{remark}

{$\bullet$} Subcase 2c: $L_2=L_{\max}$. By applying \eqref{HL_L12}, we get
\begin{align*}
\left|J(\hat{u}_{1}, \hat{u}_{2}, \hat{v})\right|
\lesssim & L_{1}^{\frac{1}{r}}L^{\frac{1}{r^\prime}}N_{2}^{\frac{1}{r}-\frac{1}{r^{\prime}}}\left(N_{1}^{\alpha}N_{2}\right)^{-\frac{1}{r^{\prime}}}\|\hat{u}_{1}\|_{L^{r^{\prime}}} \|\hat{u}_{2}\|_{L^{r^{\prime}}} \|\hat{v}\|_{L^{r}} \\
\lesssim &L_{1}^{\frac{1}{r}}L^{\frac{1}{r^\prime}}N_{2}^{\frac{1}{r}-\frac{1}{r^{\prime}}}\left(N_{1}^{\alpha}N_{2}\right)^{-\frac{1}{r^{\prime}}}L_{1}^{-\frac{1}{r}-}N_{1}^{-s}(N_1^{1+\alpha}N_{2})^{-\frac{1}{r}}N_{2}^{-s}L^{-\frac{1}{r^{\prime}}-}N_1^{+}N^{s+1}\\
&\cdot\|u_{1}\|_{\X_r^{s, 1/r+}}\|u_{2}\|_{\X_r^{s, 1/r+}}\|v\|_{\X_{r^{\prime}}^{-s-1, 1/r^{\prime}--}}\\
\lesssim & N_{1}^{1-\alpha-\frac{1}{r}+}N_{2}^{-s-\frac{2}{r^\prime}}\|u_{1}\|_{\X_r^{s, 1/r+}}\|u_{2}\|_{\X_r^{s, 1/r+}}\|v\|_{\X_{r^{\prime}}^{-s-1, 1/r^{\prime}--}}.
\end{align*}
We find the condition $r<\frac{1}{1-\alpha}$, which is weaker than $r<1+\alpha$ for $0<\alpha \leq 1$. Hence, we demand $1<r < 1+\alpha$ and $s >-1-\alpha+\frac{1}{r}$.
 
{$\bullet$} Subcase 2d: $L_{1}=L_{\max}$. Similarly, applying \eqref{HL_L12} gives
\begin{align*}
\left|J(\hat{u}_{1}, \hat{u}_{2}, \hat{v})\right|
\lesssim & L_{2}^{\frac{1}{r}}L^{\frac{1}{r^\prime}}N_{2}^{\frac{1}{r}-\frac{1}{r^{\prime}}}N_{1}^{-\frac{1+\alpha}{r^{\prime}}}\|\hat{u}_{1}\|_{L^{r^{\prime}}} \|\hat{u}_{2}\|_{L^{r^{\prime}}}  \|\hat{v}\|_{L^{r}} \\
\lesssim &L_{2}^{\frac{1}{r}}L^{\frac{1}{r^\prime}}N_{2}^{\frac{1}{r}-\frac{1}{r^{\prime}}}N_{1}^{-\frac{\alpha+1}{r^{\prime}}}(N_1^{1+\alpha}N_{2})^{-\frac{1}{r}}N_{1}^{-s}L_{2}^{-\frac{1}{r}-}N_{2}^{-s}L^{-\frac{1}{r^{\prime}}-}N_{1}^{+}N^{s+1}\\
&\cdot\|u_{1}\|_{\X_r^{s, 1/r+}}\|u_{2}\|_{\X_r^{s, 1/r+}}\|v\|_{\X_{r^{\prime}}^{-s-1, 1/r^{\prime}--}}\\
\lesssim & N_{1}^{-\alpha+}N_{2}^{-s-\frac{1}{r^\prime}}\|u_{1}\|_{\X_r^{s, 1/r+}}\|u_{2}\|_{\X_r^{s, 1/r+}}\|v\|_{\X_{r^{\prime}}^{-s-1, 1/r^{\prime}--}},
\end{align*}
We end up again with the condition $\alpha >0$ and  $s > -1-\alpha+\frac{1}{r}$.

As a consequence, given $0<\alpha \leq 1$, the $high \times low$ interactions give the conditions $1 <r < 1+\alpha$ and $s >-1-\alpha+\frac{1}{r}$.

{\bf Case 3:} $N_{1}\sim N_{2}\gg N$. It implies that $\xi_1$ and $\xi_2$ have different signs $(\xi_1\xi_2<0)$.
In this case, we have $|\Omega\left(\xi_{1}, \xi_{2}\right)|\ges |\xi_1|^{1+\alpha}|\xi_1+\xi_{2}|\sim N_1^{1+\alpha}N$. We begin with $N \lesssim 1$. Similarly, abandoning the use of $L_{\max}\ges|\Omega\left(\xi_{1}, \xi_{2}\right)|$ gives a better regularity for $N \lesssim 1$. Hence,

{$\bullet$} Subcase 3a: $N\lesssim 1$. Using \eqref{HH_L12} gives
\begin{align*}
\left|J(\hat{u}_{1}, \hat{u}_{2}, \hat{v})\right|
\lesssim & L_{1}^{\frac{1}{r}}L^{\frac{1}{r^\prime}}N_{1}^{\frac{1}{r}-\frac{1}{r^{\prime}}}N_{1}^{-\frac{1+\alpha}{r^\prime}}\|\hat{u}_{1}\|_{L^{r^{\prime}}} \|\hat{u}_{2}\|_{L^{r^{\prime}}}  \|\hat{v}\|_{L^{r}} \\
\lesssim &L_{1}^{\frac{1}{r}}L^{\frac{1}{r^\prime}}N_{1}^{\frac{1}{r}-\frac{1}{r^{\prime}}}N_{1}^{-\frac{1+\alpha}{r^\prime}}L_{1}^{-\frac{1}{r}-}N_{1}^{-s}L_2^{-\frac{1}{r}-}N_{2}^{-s}L^{-\frac{1}{r^{\prime}}-}N_{1}^{+}\\
 &\cdot\|u_{1}\|_{\X_r^{s, 1/r+}}\|u_{2}\|_{\X_r^{s, 1/r+}}\|v\|_{\X_{r^{\prime}}^{-s-1, 1/r^{\prime}--}}\\
\lesssim & N_{1}^{-2s-\alpha-2+\frac{2}{r}+}\|u_{1}\|_{\X_r^{s, 1/r+}}\|u_{2}\|_{\X_r^{s, 1/r+}}\|v\|_{\X_{r^{\prime}}^{-s-1, 1/r^{\prime}--}}.
\end{align*}
The boundedness is fine if $s> -1-\frac{\alpha}{2}+\frac{1}{r}$. Then it remains to consider $N \gg 1$.

{$\bullet$} Subcase 3b: $L=L_{\max}$. Using \eqref{HHL} gives
\begin{align*}
\left|J(\hat{u}_{1}, \hat{u}_{2}, \hat{v})\right|
\lesssim & L_{1}^{\frac{1}{r}}L_{2}^{\frac{1}{r}}(N_{1}^{\alpha}N)^{-\frac{1}{r}}\|\hat{u}_{1}\|_{L^{r^{\prime}}} \|\hat{u}_{2}\|_{L^{r^{\prime}}}  \|\hat{v}\|_{L^{r}} \\
\lesssim &L_{1}^{\frac{1}{r}}L_{2}^{\frac{1}{r}}(N_{1}^{\alpha}N)^{-\frac{1}{r}}L_{1}^{-\frac{1}{r}-}N_{1}^{-s}L_{2}^{-\frac{1}{r}-}N_{2}^{-s}L^{-\frac{1}{r^{\prime}}-}N_{1}^{+}N^{s+1}\\
&\cdot\|u_{1}\|_{\X_r^{s, 1/r+}}\|u_{2}\|_{\X_r^{s, 1/r+}}\|v\|_{\X_{r^{\prime}}^{-s-1, 1/r^{\prime}--}}\\
\lesssim & L_{1}^{\frac{1}{r}}L_{2}^{\frac{1}{r}}(N_{1}^{\alpha}N)^{-\frac{1}{r}}L_{1}^{-\frac{1}{r}-}N_{1}^{-s}L_{2}^{-\frac{1}{r}-}N_{2}^{-s}(N_1^{1+\alpha}N)^{-\frac{1}{r^{\prime}}}N_{1}^{+}N^{s+1}\\
&\cdot\|u_{1}\|_{\X_r^{s, 1/r+}}\|u_{2}\|_{\X_r^{s, 1/r+}}\|v\|_{\X_{r^{\prime}}^{-s-1, 1/r^{\prime}--}}\\
\lesssim& N_{1}^{-2s-\alpha-1+\frac{1}{r}+}N^{s
}\prod^2_{i=1}\|u_{i}\|_{\X_{r}^{s, 1/r+}}\|v\|_{\X_{r^{\prime}}^{-s-1, 1/r^{\prime}--}}.
\end{align*}
We find the condition $s > \frac{1}{2}\left(-1-\alpha+\frac{1}{r}\right)$. 

{$\bullet$} Subcase 3c: $L_{2}=L_{\max}$. Using \eqref{HH_L12} gives
\begin{align*}
\left|J(\hat{u}_{1}, \hat{u}_{2}, \hat{v})\right|
\lesssim & L_{1}^{\frac{1}{r}}L^{\frac{1}{r^\prime}}N_{1}^{\frac{1}{r}-\frac{1}{r^{\prime}}}N_{1}^{-\frac{1+\alpha}{r^\prime}}\|\hat{u}_{1}\|_{L^{r^{\prime}}} \|\hat{u}_{2}\|_{L^{r^{\prime}}}  \|\hat{v}\|_{L^{r}} \\
\lesssim &L_{1}^{\frac{1}{r}}L^{\frac{1}{r^\prime}}N_{1}^{\frac{1}{r}-\frac{1}{r^{\prime}}}N_{1}^{-\frac{1+\alpha}{r^\prime}}L_{1}^{-\frac{1}{r}-}N_{1}^{-s}(N_1^{1+\alpha}N)^{-\frac{1}{r}}N_{2}^{-s}L^{-\frac{1}{r^{\prime}}-}N_{1}^{+}N^{s+1}\\
 &\cdot\|u_{1}\|_{\X_r^{s, 1/r+}}\|u_{2}\|_{\X_r^{s, 1/r+}}\|v\|_{\X_{r^{\prime}}^{-s-1, 1/r^{\prime}--}}\\
\lesssim & N_{1}^{-2s-\alpha-2+\frac{2}{r}+}N^{s+1-\frac{1}{r}}\|u_{1}\|_{\X_r^{s, 1/r+}}\|u_{2}\|_{\X_r^{s, 1/r+}}\|v\|_{\X_{r^{\prime}}^{-s-1, 1/r^{\prime}--}},
\end{align*}
which is fine if $s> -1-\frac{\alpha}{2}+\frac{1}{r}$. 

{$\bullet$} Subcase 3d: $L_{1}=L_{\max}$. Since $\xi_1$ and $\xi_2$ are symmetric here, the estimate follows the same lines as in subcase 3c and leads again to condition $s> -1-\frac{\alpha}{2}+\frac{1}{r}$. 

In conclusion, the $high \times high\rightarrow low$ interactions deduce that the boundedness is fine provided that $s > \frac{1}{2}\left(-1-\alpha+\frac{1}{r}\right)$.

{\bf Case 4:} $N_{1}\sim N_{2}\sim N \gg 1$. 
In this case, we have $|\Omega\left(\xi_{1}, \xi_{2}\right)|\ges N_1^{2+\alpha}$. 

{$\bullet$} Subcase 4a: $L=L_{\max}$. The condition here relies heavily on the assumption $\xi_1 \xi_2 >0$, instead of $\xi_1 \xi_2 <0$. Hence, by applying \eqref{HHH}, we have
\begin{align*}
\left|J(\hat{u}_{1}, \hat{u}_{2}, \hat{v})\right|
\lesssim & N_{1}^{-\frac{\alpha}{2r}}\|{u}_{1}\|_{\X_{r}^{0, 1/r+}} \|{u}_{2}\|_{\X_{r}^{0, 1/r+}} \|\hat{v}\|_{L^{r}} \\
\lesssim & N_{1}^{-\frac{\alpha}{2r}}N_{1}^{-s}N_{2}^{-s}(N_1^{2+\alpha})^{-\frac{1}{r^{\prime}}}N_{1}^{+}N^{s+1}\\
&\cdot\|u_{1}\|_{\X_r^{s, 1/r+}}\|u_{2}\|_{\X_r^{s, 1/r+}}\|v\|_{\X_{r^{\prime}}^{-s-1, 1/r^{\prime}--}}\\
\lesssim& N_{1}^{-s-\alpha-1+\frac{2}{r}+\frac{\alpha}{2r}+}\prod^2_{i=1}\|u_{i}\|_{\X_{r}^{s, 1/r+}}\|v\|_{\X_{r^{\prime}}^{-s-1, 1/r^{\prime}--}},
\end{align*}
We find the condition $s> -1-\alpha+\frac{2}{r}+\frac{\alpha}{2r}$.

{$\bullet$} Subcase 4b: $L_{1}=L_{\max}$ or $L_{2}=L_{\max}$. These two cases follow a similar argument as in Subcase 3c, and lead to the condition $s> -1-\alpha+\frac{1}{r}$. They are the more harmless cases since the gain from the resonance relation lies completely on the high frequency. Therefore, we omit the proof. Finally, we end up with the demand $s> -1-\alpha+\frac{2}{r}+\frac{\alpha}{2r}$ in $high \times high\rightarrow high$ interactions.

Compiling all cases above and using the standard contraction mapping argument give the desired result of Theorem \ref{theorem}.
\end{proof}

\section*{Acknowledgements}
I am grateful for the support and patience of my supervisor, Zihua Guo, who provided useful discussions throughout this paper. Moreover, I would like to thank Professor Axel Gr\"unrock. His papers greatly inspired me and he provided me with valuable guidance for improving the key estimates. In addition, I would like to thank Professor Sebastian Herr for useful comments.

\end{document}